%
%
\magnification=1000

\baselineskip=14.1pt

\hsize = 15truecm \vsize = 22truecm
\hoffset = 0.4truecm \voffset = 0.7truecm

%
%

\font\smallcaps=cmcsc10

\font\title=cmr10 scaled 1200

%
%
\def\qed{${\vcenter{\vbox{\hrule height .4pt
           \hbox{\vrule width .4pt height 4pt
            \kern 4pt \vrule width .4pt}
             \hrule height .4pt}}}$}
\def\mqed{{\vcenter{\vbox{\hrule height .4pt
           \hbox{\vrule width .4pt height 4pt
            \kern 4pt \vrule width .4pt}
             \hrule height .4pt}}}}
%
%
\font\bbb=msbm10 at 10pt  

\def\RR{\hbox{\bbb R}}
\def\CC{\hbox{\bbb C}}
\def\ZZ{\hbox{\bbb Z}}
\def\PP{\hbox{\bbb P}}
\def\QQ{\hbox{\bbb Q}}

\def\hash{\#}

\newcount\refCount
\def\newref#1 {\advance\refCount by 1
\expandafter\edef\csname#1\endcsname{\the\refCount}}
\newcount\derandomCount
\def\derandom#1 {\advance\derandomCount by 1
\expandafter\edef\csname#1\endcsname{\the\derandomCount}}
%
%

%
%

\centerline{\title
Rational points in periodic analytic sets and the Manin-Mumford
conjecture}

\bigskip

\centerline{Jonathan Pila \quad \& \quad Umberto Zannier}\bigskip\bigskip

\noindent{\smallcaps Abstract.} We present a new proof of the Manin-Mumford conjecture about torsion points on algebraic subvarieties of abelian varieties. Our principle, which admits other applications,  is to view torsion points as rational points on a complex torus and then compare   (i) upper bounds for the number of rational points on a  transcendental analytic variety (Bombieri-Pila-Wilkie) and (ii) lower bounds for the degree of a torsion point (Masser), after taking conjugates. In order to be able to deal with  (i), we  discuss (Thm. 2.1) the semi-algebraic curves contained in an  analytic variety supposed  invariant for translations by a full lattice, which is  a topic with some independent motivation.
\bigskip

\noindent{\smallcaps \S 1. Introduction}\medskip 

The so-called {\it Manin-Mumford
conjecture} was raised independently by Manin and Mumford and first proved by
Raynaud [R] in 1983; its original form stated that {\it a curve $C$ (over $\CC$) 
of
genus $\ge 2$, embedded in its Jacobian $J$, can contain only finitely many 
torsion
points} (relative of course to the Jacobian group-structure). Raynaud actually
considered the more general case when
$C$ is embedded in any abelian variety. Soon afterwards, Raynaud [R2] produced a 
further significant generalization,
replacing $C$ and $J$  respectively by a subvariety $X$ in an abelian 
variety $A$; in this situation
he proved  that {\it if $X$ contains a Zariski-dense set of torsion points, then
$X$ is a translate of an abelian subvariety
of $A$ by a torsion point}. Other proofs (sometimes only for the case of curves) later appeared, due to Serre, to  
Coleman, to Hindry, to Buium, to Hrushovski (see [Py]),  to Pink \& Roessler [PR], to M. Baker \& Ribet [BR]. We also remark that a less deep
precedent of this problem was an analogous question for  multiplicative 
algebraic groups, raised by Lang already in
the '60s. (See [L]; Lang started the  matter by asking to describe  
the plane curves $f(x,y)=0$ with
infinitely many points $(\zeta,\eta)$, with $\zeta,\eta$ roots of unity.) 

In the meantime, the statement was put into a
broader perspective, by viewing it as a special case of the general Mordell-Lang 
conjecture and also, from another
viewpoint, of the Bogomolov conjecture on points of small canonical height on 
(semi)abelian varieties (we recall that
torsion points are those of zero height).
 These conjectures have later been proved and unified (by Faltings, Vojta,
Ullmo, Szpiro, Zhang, Poonen, David, Philippon,...) 
by means of different approaches providing, as a byproduct,
several further proofs of the Manin-Mumford statement
(we refer to the survey papers [Py] and [T]  for  a history of the topic and for
references.) 
Recent work of Klingler, Ullmo and Yafaev proving (under GRH) the Andr\'e-Oort conjecture, an analogue of the Manin-Mumford
conjecture for Shimura varieties, has inspired another proof of Manin-Mumford
due to Ratazzi \& Ullmo [RU]. All of these approaches are rather sophisticated and
depend on tools of various nature. 

\medskip 

It is the purpose of this paper to present a 
completely different proof compared to the existing
ones. Our approach too relies on certain auxiliary  results, having  
however   another  nature (archimedean)  with respect  to
the prerequisites of the previously known proofs; hence we believe that this
treatment  may be of some interest for a number
of mathematicians. Also, and perhaps more importantly, the underlying 
principle has certainly other applications, as
in work in progress [MZ]; we shall say a little more on this at the end. 

In short, the basic strategy of our proof is
as follows: view the torsion points as rational points on a real torus; 
estimate from above the number of rational
points on a transcendental subvariety (Pila-Wilkie); estimate from below the
number of torsion
 points by considering degree and
taking conjugates (Masser); obtain a contradiction if the order of torsion is 
large. In more detail, we can proceed
along the following steps, sticking for simplicity to the case of an algebraic 
curve $X$, of genus $\ge 2$, in the
abelian variety $A$, both over a number field. (The general case of complex
 numbers can be dealt with in a similar way
or reduced to this by specialization.) 

\bigbreak (i) There is a complex analytic group-isomorphism $\beta:
\CC^g/\Lambda\to A$, where $\Lambda$ is a certain lattice of rank $2g$, 
say with a $\ZZ$-basis $\lambda_1,\ldots
,\lambda_{2g}$; so we can view a torsion point $P\in A$ as the image $P=\beta(x)$ where $x$ is the class modulo
$\Lambda$ of a vector $r_1\lambda_1+\ldots +r_{2g}\lambda_{2g}\in\CC^g$, 
where the $r_i$ are rationals; if $P$ has
exact order $T$, the $r_i$ will have exact common denominator $T$. 

(ii) The algebraic curve $X\subset A$ equals
$\beta(Y)$, where $Y=\beta^{-1}(X)$ is a complex analytic curve in
$\CC^g/\Lambda$; in turn,  if $\pi:\CC^g\to \CC^g/\Lambda$ is the
natural projection, we may write $Y=\pi(Z)$, where $Z=\pi^{-1}(Y)\subset\CC^g$ is
an  analytic curve which is invariant under
translations in $\Lambda$.

 (iii) We can use the basis 
$\lambda_1,\ldots,\lambda_{2g}$ to view $\CC^g$ as $\RR^{2g}$;
then $Z$ will become a real analytic surface
in $\RR^{2g}$, denoted again $Z$, invariant under $\ZZ^{2g}$. Also, as in
(i), the torsion points on $A$ will correspond to rational points in $\RR^{2g}$. 
Then, the torsion points on $X$ will
correspond to rational points on $Z$. Note that, in view of the invariance of 
$Z$ by integral translations, it
suffices to study the rational points on $Z$ in a bounded region of $\RR^{2g}$. 

(iv) Due to a method introduced by
Bombieri-Pila [BP] for curves, and further developed by Pila [P] 
for surfaces, and by Pila-Wilkie [PW] in higher
dimensions, one can get good estimates for the number of rational points 
with denominator $T$ on a bounded region of a
real analytic variety. As has to be expected, 
these estimates apply only if one confines the attention to the rational
points which do not lie on any of the real semi-algebraic curves on the variety. 
For the number of these points, the
estimates take the shape $\ll T^\epsilon$, for any given $\epsilon >0$. This 
can be applied to the above defined $Z$;
it turns out that, since $X$ is not a translate of an elliptic curve, 
$Z$ does not contain any semi-algebraic curve.

(v) All of this still does not still yield any finiteness result, 
merely estimates for the number of torsion points on
$X$, of a given order $T$. The crucial issue is that if $X$ contains an
algebraic   point $P$, it automatically contains its
conjugates over a field of definition. Now, Masser has proved in 1984 
that the degree, over a field of definition for
$A$, of any torsion point of exact order $T$ is 
$\gg T^\rho$ for a certain $\rho>0$ depending only on $\dim A$. Hence,
if $X$ contains a point of order $T$, it contains at least 
$\gg T^\rho$ such points. 

Then, comparing the estimates
coming from (iv) and (v) we deduce that the order of the torsion points on 
$X$ is bounded, concluding the argument.\medskip

In some previous proofs, one exploited not a lower bound for degrees, but rather the Galois structure of the field generated by torsion points. However this may be considered an information substantially of different nature. 

To better isolate the new arguments from previous ones, we will prove the
Manin-Mumford conjecture in the following weak form, which includes the case of
curves.
\medskip 
\noindent {\bf  Theorem 1.1.\/} {\it Let $A$ be an abelian variety and
$X$  an algebraic subvariety of $A$, both defined over a numberfield. Suppose
that
$X$ does not contain any translate of an abelian subvariety of $A$ of dimension
$>0$. Then $X$ contains only finitely many torsion points of $A$.\/} \medskip 

However, 
rather elementary purely geometrical considerations based on degrees allow one
to get easily the following more precise version, for which we shall give a
sketch at the end; in this statement by `torsion coset' we mean a translate of an abelian subvariety by a torsion point: 
\medskip 

\noindent {\bf \ Theorem  1.1*.\/} {\it Let $A$
 be an abelian variety and $X$ an algebraic subvariety of $A$, both 
defined over a numberfield. The Zariski closure of the set of
torsion points of $A$ contained in $X$ is a finite union of torsion cosets.\/}
 \medskip 

Concerning the proof-method, a further point is that, in order for the Pila-Wilkie estimates to be applicable, we must study
the real semi-algebraic curves that may lie on the inverse image $Z$ of $X$
under the analytic uniformization $\CC^g\rightarrow A$. The set $Z$ is analytic
(defined by the vanishing of certain polynomials in the abelian functions giving
the map $\CC^g\rightarrow A$) and periodic modulo the lattice $\Lambda$. 

In \S2,
Theorem 2.1, we show that a connected real semi-algebraic curve contained in
such a set $Z$ must be contained in a complex linear subspace contained in $Z$,
and moreover a subspace in which the period lattice has full rank. In other
words, we prove that the ``algebraic part'' of $Z$ in the sense of Pila-Wilkie
(see Def. 2.1 below) corresponds precisely to the union of translates of abelian
subvarieties of $A$ of dimension $>0$ and 
 contained in $X$. This fact seems to be not entirely free of independent
 interest. 
\medskip

\noindent{\bf Remark.} One can also check that our proof of Theorem 2.1 works in
fact for arbitrary complex tori $\CC^g/\Lambda$, even if they are not projective.
\medskip

With this result in hand the proof
of Theorem 1.1 is concluded in \S3 by combining the ingredients mentioned above.

\bigskip 
\noindent {\smallcaps \S 2. Structure of the algebraic part of a
periodic analytic set}
 \medskip 

 As announced, this section will be
devoted to a complete description of algebraic subsets of an analytic subvariety
of a complex torus, which corresponds to an analytic subvariety of $\CC^g$
periodic for a full lattice $\Lambda$, that is, invariant by translations in
$\Lambda$. We express explicitly all of this in a few definitions and then 
state the main Theorem 2.1. 
\medskip

We take a (full rank) lattice $\Lambda\subset\CC^g$ which 
will be fixed throughout. We further fix a $\ZZ$-basis
$\lambda_1,\ldots,\lambda_{2g}$ of $\Lambda$ and use it to 
identify $\CC^g$ with $\RR^{2g}$. 

We shall use throughout
the notation $Z+\Lambda:=\cup_{\lambda\in\Lambda} 
(Z+\lambda ):=\{z+\lambda\ :\ z\in Z, \lambda\in\Lambda\}$ for the
union of the translates of $Z$ by the vectors in $\Lambda$. A set 
$Z\subset\CC^g$ will be called {\it ($\Lambda$-)
periodic\/} if $Z+\Lambda=Z$. 

As usual, an {\it analytic\/} set $Z\subset \CC^g$ will mean a set such that
every point
$z\in \CC^g$ has an open neighbourhood $U$ in which $Z$ is defined as the set of 
common zeros of a finite collection of
functions (depending on $z$) that are (complex) analytic (i.e. regular) in $U$. Such a 
set is readily seen to be {\it real
analytic\/} as a subset of $\RR^{2g}$. 
\medskip

 \noindent {\bf Definition 2.1.\/} Let $Z\subset\CC^g$. We define the
{\it complex algebraic part\/} $Z^{\rm ca}$ of $Z$ to be the union of 
connected closed algebraic subsets of $\CC^g$ of
positive dimension contained in $Z$. Viewing $Z$ as a subset of 
$\RR^{2g}$ we define the {\it real algebraic part\/}
$Z^{\rm ra}$ of $Z$ to be the union of connected real algebraic sets 
of positive dimension contained in $Z$. Finally
we define the {\it algebraic part\/} $Z^{\rm alg}$ of $Z$ to be the 
union of all connected real semi-algebraic sets of
positive dimension (see [Sh, pp. 51, 100] for definitions) contained in $Z$. 
One readily sees that $Z^{\rm ca}\subset
Z^{\rm ra}\subset Z^{\rm alg}$. 
\medskip 

We reserve the term {\it subspace\/} for a (respectively complex or
real) vector subspace of $\CC^g$ or $\RR^{2g}$. By a (respectively complex or
real) {\it linear subvariety\/} we mean a subset of $\CC^g$ or $\RR^{2g}$
defined by the vanishing of some (complex or real) linear equations, not
necessarily homogeneous. A subspace $H$ in which $H\cap\Lambda$ has full rank in
$H$ will be called a {\it ($\Lambda$-) full\/} subspace. A set of the form $z+H$
where $z\in\CC^g$ and $H$ is a complex linear subspace will be called a {\it
coset\/}. Thus a subspace $H$ that is both complex and full corresponds
precisely to a {\it subtorus\/} of
$\CC^g/\Lambda$, and a coset of such an $H$ will be called a {\it torus
coset\/}.\medskip

\noindent {\bf Definition 2.2.\/} For a set $Z\subset\CC^g$ we let $Z^{\rm torus\
coset}$ be the union of torus cosets of positive dimension contained in $Z$. 
\medskip 

A torus coset is evidently a complex linear subvariety, whence $Z^{\rm torus\
coset}\subset Z^{\rm ca}\subset Z^{\rm ra}\subset Z^{\rm alg}$. For a periodic
analytic set
$Z$ we show that all these sets coincide. 
\medskip 

\noindent {\bf Theorem 2.1.\/} {\it Let $Z\subset\CC^g$ be a periodic analytic
set. Then $Z^{\rm alg}=Z^{\rm torus\ coset}$.\/} 
\medskip 

Our proof of this
result involves several preliminaries; for the reader's convenience we
explicitly subdivide the proof into three Steps. 
\medskip 

{\bf Step 1.} {\it Reduction semi-algebraic
$\to$ complex algebraic.} In this step we reduce to complex 
curves contained in $Z$, by proving that $Z^{\rm
alg}=Z^{\rm ca}$. We first prove a lemma:
\medskip

\noindent {\bf Lemma 2.1.\/} {\it Let $Z\subset\CC^g$ be analytic.
Suppose that $x\in\CC^g$ has a neighbourhood $U$ such that $x$ is a 
smooth point of $Y\cap U$, where $Y$ is a real
algebraic curve with $Y\cap U\subset Z$ as subsets of $\RR^{2g}$. Then 
there is a neighbourhood $U'$ of $x$ contained
in $U$ and a finite collection of 
irreducible complex algebraic curves $\Gamma_j$ such that
$Y\cap U'\subset\cup_j\Gamma_j\subset Z$.\/} 
\medskip 

\noindent {\it Proof.\/} We took coordinates in $\RR^{2g}$ using the
lattice basis, but clearly
$Y$ will remain semi-algebraic under any real 
linear change of coordinates. Also, for our purposes we can assume that $Y$ is
irreducible. 

Let us here write
$(x_1, y_1,\ldots,x_n, y_n)$ for the coordinates of $\RR^{2g}$, where
$z_j=x_j+iy_j$  are the coordinates of $\CC^g$. We may assume by
translation that $x=0$. If $x_1$ (say) is a non-constant function on the real 
curve $Y$ near 0 then, for $t=x_1$ in
some real neighbourhood of 0, all the functions $x_1, y_1,\ldots,x_n, y_n$
are  real analytic functions of some $m$-th root
$u:=t^{1/m}$, and (as functions of $u$), are algebraic over 
$\RR(x_1)\subset\CC(x_1, y_1)$. We claim that each of the
functions $x_2, y_2,\ldots,x_n, y_n$ is algebraic over the 
field $\CC(z_1)=\CC(x_1+iy_1)$. This is because the
function $x_1+iy_1$ is non-constant on $Y$ since $x_1$ is non-constant, 
and $\CC(x_1, y_1)$ is algebraic over
$\CC(x_1)$, so of transcendence degree 1. Therefore, for each $j$, 
the function $z_j=x_j+iy_j$ is also algebraic over
$\CC(z_1)$. The functions $z_j(u)$ are analytic for real $u$ in some
 neighbourhood of 0, hence for complex $u$ in a
complex neighbourhood of 0;  the image of 
 $u\mapsto z(u)$ is thus a complex open subset of a complex algebraic
irreducible  curve $\Gamma$, which in a neighborhood of $0$ contains $Y$ and 
must
necessarily be contained in
$Z$.  Certainly $\Gamma$ contains a smooth point in this neighborhood, and we get
$\Gamma\subset Z$ by analytic continuation, as the smooth points of an 
irreducible complex curve are connected, and the remaining points
are isolated and belong to $Z$ by continuity. 
\ \qed 
\medskip

We can now prove the alluded
reduction:
\medskip

 \noindent {\bf Proposition 2.1.\/} {\it Let $Z\subset\CC^g$ be a
periodic analytic set. Then
$Z^{\rm alg}=Z^{\rm ca}$.}
 \medskip 

\noindent{\it Proof.} 
Suppose $W$ is a connected real semi-algebraic set of
positive dimension with $W\subset Z$. Then,   
omitting at most finitely many points, $W$ is a union of real connected
semi-algebraic curves. If $Y$ is such a curve then,  $Y$ is real algebraic
in the neighbourhood  of any smooth point and   we may apply  Lemma 2.1 to
find $Y$ contained in a complex algebraic
 curve $\Gamma\subset Z$. This proves the
statement. \ \qed 
\medskip

{\bf Step 2.} {\it Complex curves in periodic analytic varieties and
linear spaces.} In this step we prove that if a complex curve $C$ is contained
in the periodic analytic set $Z$, then there are `several' complex lines $l$
with $C+l$ contained in $Z$. In turn, this will involve a few preliminaries.
\medskip Suppose
$C\subset\CC^g$ is an irreducible complex algebraic curve. Of the
 coordinates $z_1,\ldots,z_g$ of $\CC^g$, if $z_i$ is  not constant on $C$
then any $z_i,z_j$  are  related by some irreducible polynomial equation 
$G(z_i, z_j)=0$. For $|z_i|$ sufficiently large, say $|z_i|>R$, the 
solutions $z_j$ are given by convergent Puiseux
series $\phi_j(z_i)$. By a {\it branch\/} of $C$ we mean a choice of index $i$
and a 
$g$-tuple $\phi=(\phi_1(z_i),\ldots,
\phi_g(z_i))$ of algebraic Puiseux series, convergent for $|z_i|>R$, and 
such that
$\phi(z_i)=(\phi_1(z_i),\ldots,\phi_n(z_i))\in C$ for all $|z_i|> R$. By a 
suitable choice of index $i$, we can always
obtain a branch of $C$ such that, for all $j$, 
$$ 
\phi_j(z_i)=\alpha_jz_i+{\rm\ lower\ order\ terms},\qquad
\alpha_j\in\CC. 
$$
 We call such a branch {\it linear\/}, and 
$\alpha=(\alpha_1,\ldots,\alpha_g)$ the {\it direction\/}
of the branch. Observe that $\alpha_i=1$, so $\alpha\ne 0$.

Suppose $\phi(w)=(\phi_1(w),\ldots,\phi_g(w))$ is an
$n$-tuple of algebraic Puiseux series, convergent for $|w|>R$. 
Then, fixing $w_0, w_1\in\CC$ and
$\mu=(\mu_1,\ldots,\mu_g)\in\CC^g$, we consider (for $\kappa\in\CC$ 
such that $|w_0+\kappa w_1|>R$), 
$$ 
\kappa\mapsto
\psi(\kappa)=(\psi_1(\kappa),\ldots,\psi_g(\kappa)), 
$$ 
$$ 
\psi_i(\kappa)=\phi_i(w_0+\kappa w_1)-\kappa\mu_i. 
$$
 From the  algebraic relation $G_i(t,\phi_i(t))$, we
 find that $G_i(w_0+\kappa w_1,
\psi_i(\kappa)+\kappa\mu_i)=0$. Thus the $\psi_i(\kappa)$ are also 
algebraic functions of $\kappa$, and the locus above is Zariski dense in an
irreducible algebraic curve in $\CC^g$, which we denote $\Gamma(\phi,\mu, w_0,
w_1)$. It is important to note that its degree is bounded in terms of the 
degree of the curve containing $\phi$ (independently of the choice of 
$w_0,w_1, \mu$). We now have a lemma, perhaps
known but for which we have found no reference: 
\medskip 

\noindent {\bf Lemma 2.2.\/} {\it Let $Z\subset \CC^g$ be an
analytic set, $B\subset\CC^g$ a bounded set, and $\delta$ a positive integer. 
There exists a positive integer
$K=K(Z,B,\delta)$ with the following property. Suppose $\Gamma\subset\CC^g$ 
is an irreducible complex algebraic curve
of degree $\le \delta$, and with $\hash (Z\cap B\cap \Gamma)\ge K$. Then 
$\Gamma\subset Z$.} 
\medskip 

\noindent {\it Proof.\/} About each point of $\CC^g$  there is an open disk in which
$Z$ is  defined by the vanishing of a finite number of regular functions. Taking a smaller disk,
these functions may be assumed regular in a neighbourhood of the closure of the disk. 
Then by compactness of the closure of $B$, we can find a finite number of open disks $U$ covering $B$, and  on each disk a finite number of  functions, regular in a neighbourhood of the closure of the disk, 
whose zero-locus defines $Z$ in the disk. Thus, in each of the
finitely many disks, $Z$ is expressed as in intersection of finitely many 
hypersurfaces. As shown in a moment, this remark  reduces the proof to the case
of hypersurfaces, namely to the following
 \medskip

 \noindent {\bf Claim.\/} {\it Suppose that $U\subset\CC^n$ is a bounded open 
disk and that $f$ is a complex-valued function that is regular in a neighbourhood of the closure of $U$. 
Let $Y=\{z\in U: f(z)=0\}$, and $\delta$ a positive integer.
There is a positive integer $K=K(U, f, \delta)$ with the following property.
 Suppose $\Gamma\subset\CC^n$ is an
irreducible curve of degree $\le \delta$ and $\hash Y\cap C\ge K$. 
Then $f$ vanishes identically on $\Gamma$ in $U$.}
\medskip 

Given the claim, we may establish the lemma as follows. 
For each disk $U$ in the finite covering take $K_U$
 to be the maximum of the $K(U,f,\delta)$ over all 
the finitely many functions $f$ defining $Z$ on $U$, and take
$K=\sum_UK_U+1$. For each $U$ put $Z_U=Z\cap U$. Let now $\Gamma$ be
 an irreducible complex algebraic curve of degree $\le\delta$ with $\hash
Z\cap B\cap\Gamma\ge K$. By the pigeonhole principle, one of the 
disks $U$ has $\hash (Z_U\cap\Gamma)\ge K_U$. By the
Claim, each of the functions $f$ defining $Z$ on $U$ vanishes
 identically on $\Gamma$, so that, restricted to $U$,
$\Gamma\subset Z$. However $\Gamma$ is connected, so by analytic 
continuation we conclude that $\Gamma\subset Z$. It
remains then to prove the Claim, as we shall now do. 

\medskip

Since $\Gamma$ has degree $\le \delta$, each pair of
coordinates $z_i, z_j$ satisfy on $\Gamma$ an algebraic relation
 $P_{ij}(z_i, z_j)=0$ for some nonzero polynomial
$P_{ij}\in\CC[x_1, x_2]$ of degree $\le \delta$. These polynomials,
 considered up to nonzero constant factors, define
an algebraic set of dimension $1$, containing $\Gamma$ 
as a component. This set is  possibly reducible but, given the $P_{ij}$, by projecting to a
general plane, we see that there certainly exists a hypersurface $H$ of degree
$\le g\delta$ that contains
$\Gamma$  but not any other component of dimension $1$;
this $H$ corresponds to a further polynomial $P_H$ of degree $\le g\delta$.
 Now, the polynomials in two variables of
degree $\le \delta$, up to constants, are parametrised by a projective space
 $\PP^D$, $D=(\delta+1)(\delta+2)/2$, whereas
$H$ is parametrized by $\PP^{D'}$ for a $D'$ depending only on $g,\delta$. 
Hence the set  $\{H,P_{ij}\}$ is parametrized
by the product 
$$ 
\Delta=\PP(\CC)^{D'}\prod_{ij}\PP(\CC)^D,
$$ 
a compact space, which we can assume contained in  $\RR^m$ for some suitable $m$. 
For $w\in\Delta$ let $I_w$ be
the ideal generated by the corresponding $P_{ij}, P_H$. Consider 
the set 
$$ 
V:=\{(z,w)\in \CC^g\times\Delta:
z=(z_1,\ldots,z_g)\in U, w\in\Delta, f(z)=0, Q(z)=0, {\rm\ all\ } Q\in I_w\}.
 $$ 
We have
$V\subset\CC^g\times\Delta$
 with projections $\pi_1, \pi_2$ onto the factors $\CC^g$ and
$\Delta$ respectively.

 Our very construction shows that every algebraic 
curve $\Gamma\subset\CC^g$ of degree
$\le\delta$ is defined, up to finitely many points, by $I_w$ for (at least one) 
suitable $w=w_\Gamma\in\Delta$; hence, for some point
$w_\Gamma\in\Delta$, $\pi_1\big(\pi_2^{-1}(w_\Gamma)\big)$
equals the intersection $Y\cap \Gamma$ plus a finite set. 

Now we consider $V$ as a subset of $\RR^{M}$ for suitable $M$. 
Since $f$ is regular on a neighbourhood of the closure of $U$, the set $V$ is {\it subanalytic} (even {\it semianalytic\/})
in $\RR^M$
(for the definition see [G] or [BM, Definition 2.1]). Further, $V$ is
bounded, since $U$ is bounded and $\Delta$ is compact. We appeal to
Gabrielov's Theorem [G, or see e.g. BM, Theorem 3.14] to conclude: 
\medskip 

{\it As $w$ varies over the bounded set
 $\pi_2(V)$, the number of connected
components of $\pi_2^{-1}(w)$ is bounded by some finite number $N=N(U,f,D)$; 
the number of connected components of
$\pi_1(\pi_2^{-1}(w))$ is then also bounded by $N$. } 
\medskip 

Put $K=N+1$.
 If now $\hash Y\cap \Gamma\ge K$,
$Y\cap\Gamma$ cannot consist of isolated points and must, as a semianalytic 
set, have dimension $\ge 1$. This set must
then contain some smooth real-analytic arc. If we take a point on such an arc 
that is a non-singular point of
$\Gamma$, then restricting $f$ to a neighbourhood in which we can complex 
analytically parameterize $\Gamma$, we find
that  $f$ restricted to this local parameterization is an 
analytic function with non-isolated zeros. It
therefore vanishes identically on
$\Gamma$ in this neighbourhood. But now since $\Gamma$ is irreducible, the set
of its non-singular points is connected, and we find that $f$ vanishes
identically on
$\Gamma$ in $U$ by analytic continuation. This establishes the Claim, and
concludes the proof of Lemma 2.2. \ \qed 
\medskip 

The following result now
follows easily:
\medskip 

\noindent {\bf Corollary.\/} {\it Let $Z\subset\CC^g$ be
a periodic analytic set and
$r>0$. Let $\phi(w)$ be an $n$-tuple of convergent (for $|w|>R$) algebraic 
Puiseux series. There is an integer $K=K(Z,r,\phi)$ with the following property.
Let $B\subset\CC^g$ be a ball of radius $r$, and $w_0,
w_1\in\CC,\tau,\mu\in\CC^g$. If 
$$ 
\hash\big(\Gamma(\phi,\mu, w_0, w_1)+\tau\cap
Z\cap B\big)\ge K 
$$ then
$\Gamma(\phi,\mu, w_0, w_1)+\tau\subset Z$.}\ \qed 
\medskip 

\noindent {\it Proof.\/}
For a fixed $B$ the conclusion follows directly from Lemma 2.2, because
$\Gamma(\phi,\mu, w_0, w_1)+\tau$ is an algebraic curve of degree bounded only
in terms of $\phi$. Since
$Z$ is periodic, the ball $B$ may be assumed to be centred in a fundamental
domain, so dependence on the centre of the ball may be eliminated.\ \qed 
\medskip 

We now
come to the fundamental point of this Step 2. In the following, $||\cdot||$
denotes the Euclidean norm in $\CC^g$. For a point
$z=(z_1,\ldots,z_g)\in\CC^g-\{0\}$ we denote by $[z]$ the complex line $\{wz:
w\in\CC\}$. 
\medskip 

\noindent {\bf Proposition 2.2.\/} {\it Let $Z\subset\CC^g$
be a periodic analytic set and $C$ be an irreducible complex algebraic curve
with $C+\tau\subset Z$ for a  $\tau\in\CC^g$. Let $\phi(w)$ be a linear
branch of
$C$ (convergent for
$|w|>R$) with direction $[\alpha]$ and put $K=K(Z,1,\phi)$, as in the last
Corollary. Finally, let
$\lambda\in\Lambda$ be such that $||\beta-\lambda||<(2K)^{-1}$ for some 
$\beta\in [\alpha]$. Then $C+\tau+[\lambda]\subset Z$.\/}
\medskip
 
\noindent {\it Proof.\/} Fix $w_1\in\CC$ such that
$||w_1\alpha-\lambda||<(2K)^{-1}$. Now choose $R_1\ge R$ such that, for any $w\in\CC$
with $|w|>R_1+K|w_1|$ and any $w'\in\CC$ with $|w'|\le K|w_1|$ we have 
$$
||\phi(w+w')-\phi(w)-w'\alpha||<1/2. 
$$ 
This is possible by the condition on
$\phi$ that all terms apart perhaps from the leading term are sub-linear. If now
$|w_0|>R_1+K|w_1|$ and $k=0,1,\ldots, K$ we have 
$$
||\phi(w_0+kw_1)-\phi(w_0)-k\lambda||\le
||\phi(w_0+kw_1)-\phi(w_0)-kw_1\alpha||+k||w_1\alpha-\lambda||<1. 
$$
 Now, $\phi(w_0+kw_1)+\tau-k\lambda+\in Z$ because $Z$ is periodic and 
$C+\tau\subset Z$; also, 
$\phi(w_0+kw_1)+\tau-k\lambda\in \Gamma:=\Gamma(\phi, \lambda, w_0, w_1)+\tau$
provided
$|w_0+kw_1|>R$, and by the above all these points lie in the ball of radius 1 
about $\phi(w_0)+\tau$ for $k=0,1,\ldots,K$. 

By the Corollary to Lemma 2.2, we deduce that 
$\Gamma\subset Z$.

 We thus find $\Gamma(\phi, \lambda, w_0, w_1)+\tau\subset Z$ for
all $w_0$
 suitably large, but we must still show that in fact all $C+\tau+[\lambda]\subset
Z$.  If we set $y=w_0+\kappa w_1$, $\kappa\in\CC$, we find that 
$$ 
\phi(y)+\tau-((y-w_0)/w_1)\lambda\in Z 
$$
provided $|y|>R$. But, fixing $y,\kappa , w_1$, the above represents a line as
$w_0$ varies. Since a segment of this line lies in $Z$, the whole line does, by
analytic continuation. So for $x\in\CC$ the points 
$$ 
\phi(y)+\tau-x\lambda\in Z 
$$
provided only
$|y|>R$. Fixing now $x$, the curve in $y$ is just a branch of $C+\tau-x\lambda$.
This curve is irreducible, being a translate of $C$, and so
$C+\tau-x\lambda\subset Z$ for all $x$.\ \qed 
\medskip

{\bf Step 3.} {\it Linear subvarieties of a periodic
analytic variety.} We study general linear subvarieties of a periodic analytic
variety $Z$, their intersections with the lattice $\Lambda$, and use all of this
to exploit the important conclusion of the last proposition, which produces
certain translates $C+l$ of $C$ by a line, contained in $Z$. 
\medskip

 We   first recall  some useful simple facts from the known theory of closed
subgroups of real vector spaces, and start by noting that: {\it if
$H$ is a real subspace of $\CC^g$ then the closure of $H+\Lambda$ has the form 
$K+\Lambda$ where $K$ is a full real subspace containing $H$.}

 This statement follows from  the
description of closed subgroups of real vector spaces given in [S, Lecture VI,
\S 2], which implies that the closure in $\CC^g$ of $H+\Lambda$ has the form
$K+\Lambda_0$ for $K$ a real subspace and  $\Lambda_0$ a lattice in a space $W$
complementary to $K$. Clearly
$H\subset K$. Now, the projection of $\Lambda $ to $W$ (along $K$) is
$\Lambda_0$, which must then have full rank in
$W$, because $\RR^{2g}=\RR\Lambda\subset \RR(K+\Lambda_0)=K+\RR \Lambda_0$.
Now, lifting a basis of
$\Lambda_0$ to
$\Lambda$ we see that
$\Lambda_1:=\Lambda\cap K$ has maximal rank, i.e. equal to $\dim K$, as wanted. 

We further note that {\it for any open ball $I$ around $0$ in $K$, the set $H+I$
contains a set of generators for $\Lambda_1:=\Lambda\cap K$.}  To prove this,
let  $\Lambda'$ denote  the lattice generated by $\Lambda\cap (H+I)$ and observe
that $H+\Lambda'$ is dense in $K$: in fact, by definition of $\Lambda'$, the closure of $H+\Lambda'$ 
   contains the intersection of $I/2$ with the closure of $H+\Lambda$, so
it contains $I/2$, whence it must contain the whole $K$. Now, let
$\lambda\in\Lambda_1$; it belongs to the closure $K$ of $H+\Lambda'$; hence
$\lambda-I$ intersects $H+\Lambda'$, so
$\lambda+\Lambda'$ intersects $\Lambda\cap (H+I)\subset\Lambda'$, proving that
$\lambda\in\Lambda'$ and so $\Lambda_1\subset\Lambda'$, as desired.
\medskip

Now we may give the following
definitions:\medskip
\noindent {\bf Definition 2.3.\/} Let $H$ be a real subspace of $\CC^g$.
\smallskip 

1. We denote by ${\rm c}(H):=\CC H$ the  complex space generated by $H$ and
call it the {\it complex closure\/} of $H$;  we call $H$ {\it complex\/} if
$H={\rm c}(H)$. 

2. We denote by ${\rm f}(H)$ the {\it full closure\/} of $H$,
namely the full real subspace $K$ such that $H+\Lambda$ is dense in $K+\Lambda$,
as in the above remark. So $H$ is full just if ${\rm f}(H)=H$. 

3. We denote by
${\bf fc}(H)$ the {\it full-complex closure\/} of $Z$, namely the union of the
iterates of $H$ under the map
$H\mapsto {\rm f}({\rm c}(H))$. 
\medskip

One could alternately iterate the two operations $H\mapsto
{\rm c}(H)$, $H\mapsto {\rm f}(H)$. In general, if $g>1$ we need not have 
${\rm f}({\rm c}(H))={\rm c}({\rm f}(H))$. (In $\CC^2$ take for instance 
$H=\RR (1,0)$, $\Lambda=\ZZ (1,0)+\ZZ(i,1)+\ZZ(1,i)+\ZZ (1,\sqrt 2)$. Then
$f(H)=H$,
$c(H)=\CC (1,0)$. Since $c(H)$ is not full, we have that $fc(H)$ contains
properly $cf(H)=c(H)$.) Anyway  ${\bf fc}(H)$ may also be characterized as the
smallest subspace containing $H$ that is both full and complex. 
\medskip


\noindent {\bf Lemma
2.3.\/} {\it Let $Z\subset\CC^g$ be a periodic analytic set. Suppose $z\in\CC^g$
and $H$ is a real subspace of $\CC^g$ with $z+H\subset Z$. Then there is a torus
coset $z+M$ with $z+H\subset z+M\subset Z$, and one can take $M={\bf fc}(H)$.} 
\medskip

\noindent {\it Proof.\/} Note that by Def. 2.3(2) (which amounts to the opening
remarks of this Step 3), if
$H$ is a real subspace and $z\in \CC^g$, then $(z+H)+\Lambda$ is dense in some
$(z+K)+\Lambda$ where $K$ is a full real subspace containing $H$. Now, if $H$ is
a real subspace with $z+H\subset Z$ then, by periodicity, the fact just stated 
and continuity, we have $z+{\rm f}(H)\subset Z$, and, by analytic continuation,
we have $z+{\rm c}(H)\subset Z$. The conclusion follows by applying these
observations to $H$ and its iterates under full and complex closure.\ \qed
 \medskip


We can now rapidly conclude the proof of Theorem 2.1; we prove a last
lemma:
\medskip 

\noindent {\bf Lemma 2.4.\/} {\it Let $Z\subset\CC^g$ be a
periodic analytic set, $C$ an irreducible complex algebraic curve and $M$ a
complex subspace such that $C+M\subset Z$. Suppose that $C+M$ is not a coset of
$M$. Then there is a complex subspace $M'$ with $M\subset M'$, $\dim M<\dim M'$,
and $C+M'\subset Z$.\/} 
\medskip 

\noindent {\it Proof.\/} Let $N$ be a complex
linear subspace supplementary to $M$ in $\CC^g$. So every translate
$z+M$, where
$z\in\CC^g$, intersects $N$ in just one point.  Now, $C+M$,
 as the image of the sum-map $C\times M\to C+M$ contains an open dense
set in its Zariski closure in $\CC^g$, which is 
an irreducible algebraic
subvariety of
$\CC^g$; also, 
$(C+M)\cap N$ is irreducible (because it is the projection of $C+M$ to $N$ along
$M$). By hypothesis
$(C+M)\cap N$ is not equal to a point, hence it cannot be a finite set of
points  and therefore it contains an open dense subset of a closed irreducible
algebraic curve
$C'\subset N$; note that $C+M\subset C'+M\subset Z$ and  $\dim (C+M)=\dim M+1$. 

Take now a linear branch
$\phi$ of
$C'$ with direction
$[\alpha]$, so $[\alpha]$ lies in $N$.   Then there is  a small open ball $B$
around
$0$ in
$\CC^g$ so that $([\alpha]+B)\cap M$ does not contain any nonzero lattice point.
By a remark above, $\Lambda\cap([\alpha]+B)$ contains generators for the full
closure of
$[\alpha]$, so in particular it contains a lattice point $\lambda\not\in M$.

 Put $M'=M+[\lambda]$. 

Then $C'+M'\subset Z$. For if $\tau\in M$ we have
$C'+\tau\subset Z$, but this curve  is a translate of $C'$, and by Proposition
2.2 (applied to $C'$ in place of $C$) we find that for small enough $B$, 
$(C'+\tau)+[\lambda]\subset Z$. This being true for all $\tau\in M$, we have
$C'+M'\subset Z$. Now   so
$C+M'\subset C'+M'\subset Z$. This completes the proof.\
\qed 
\medskip

\noindent {\it Proof of Theorem 2.1\/}

 Let $C$ be  an irreducible complex algebraic curve contained
in $Z$ and take a maximal complex subspace $M$ 
 such that $C+M\subset Z$. By maximality and Lemma 2.4 we deduce that
$C+M$ is a single coset of $M$, which must be a torus coset by maximality and
Lemma 2.3. 

It follows that every closed complex algebraic set of
positive dimension contained in $Z$ is contained in the union of full cosets
of positive dimension contained in $Z$, i.e. that $Z^{\rm ca}\subset Z^{\rm
torus\ coset}$. Proposition 2.1 now finally proves what we need. \ \qed 
\medskip 

\noindent{\bf Remark.} Note that the final argument easily leads to the following
assertion: {\it If an irreducible complex algebraic set $V$ is contained in $Z$
then there is a torus coset $z+M$ with $V\subset z+M\subset Z$.} To prove this,
take a maximal torus $M$ such that for a point $z$ we have $z+M\subset Z$.
Applying the last conclusion of the proof of Theorem 2.1 to all the irreducible
curves $C$ on $V$ passing through $z$ we see by maximality that $z+M$ contains a
neighborhood of
$z$ in $V$. Hence it contains   $V$. 
\bigskip


\noindent\leftline{\smallcaps \S 3. Manin-Mumford} 
\medskip 

Let now $A$ be a
(projective) Abelian variety, and $X$ a subvariety of $A$. We have a complex
analytic uniformization $\CC^g\rightarrow A$, periodic with period lattice
$\Lambda$. The preimage $Z\subset\CC^g$ of $X$ is a periodic analytic set, and
we have shown that all real semialgebraic subsets, connected of positive
dimension, of the real reduct of $Z$ are contained in the union of torus cosets
contained in $Z$. 

Now a subtorus of $\CC^g/\Lambda$ corresponds to an abelian
subvariety of $A$ (see e.g. [R, remark on page 86], or it may be argued directly
using Chow's Theorem that such a subtorus is algebraic). 
\medskip 

\noindent {\it
Proof of Theorem 1.1.\/} Our interest is in the torsion points $P$ of $A$ that
lie on $X$. Let ${\rm tor}(A)$ denote the torsion subgroup of $A$, consisting of
all points of $A$ of finite order. The {\it order\/} $T=T(P)$ of $P$ is the
minimal positive integer with $TP=0$. A torsion point $P$ of $A$ corresponds to
a rational point  $z=z_P=(q_1,\ldots,q_{2g})\in\QQ^{2g}$ of $Z$ considered as a
subset of $\RR^{2g}$. The order of $P$ is equal to the {\it denominator\/} of
$z$, i.e. the minimal integer $d>0$ for which $dz\in\ZZ^{2g}$. By the present
assumptions,
$X$ does not contain any translate of an abelian subvariety of dimension $>0$;
hence, by Theorem 2.1 the set $Z^{\rm alg}=Z^{\rm torus\ coset}$ is empty. 

As a
subset of $\RR^{2g}$ the set $Z$ is $\ZZ^{2g}$-periodic. In considering torsion
points it therefore suffices to replace $Z$ by ${\cal Z}=Z\cap[0,1)^{2g}$, and
clearly  ${\cal Z}^{\rm alg}=Z^{\rm alg}\cap [0,1)^{2g}$ is empty as well.

For a set $W\subset [0,1)^{2g}$ and a real number $T\ge 1$ we denote by $N(W,T)$ the
number of rational points of $W$ of denominator dividing $T$. By Pila-Wilkie
[PW], for every $\epsilon>0$, 
$$ 
N({\cal Z},T)=N({\cal Z}-{\cal Z}^{\rm alg}, T)\le c_{1}({\cal
Z},\epsilon)T^\epsilon.\eqno(3.1)
$$
 On the other hand, there are lower bounds
for the degree of torsion points. Suppose $A$ is defined over a number field
$K$. For $P\in {\rm tor}(A)$ set $d(P)=[K(P):\QQ]$. Then Masser [M] proves that
$$ 
d(P)\ge c_2(A)T^{\rho} 
$$ 
for some $c_2(A)>0$ and some $\rho>0$ which depends only on the dimension $g$ of $A$. Note that all the
conjugates of $P$ over a number field of definition for both $A,X$ are still
torsion points on $X$, of the same order as $P$. By the lower bound of Masser
just displayed, the number of such conjugates is at least
$c_3(A)T^{\rho}$. Hence we get at least 
$$ 
c_3(A)T^\rho 
$$
 distinct points
$z\in {\cal Z}$ of height $\le T$, corresponding to $P$ and
its conjugates. Choosing in (3.1) $\epsilon =\rho/2$ and comparing the estimates
so obtained we conclude that $T$ is bounded.\ \qed 
\medskip 

As anticipated in the
Introduction, we note that the same arguments may be used to prove the stronger
version given by Theorem 1.1*, at the cost of adding some simple geometrical
facts on subvarieties of abelian varieties and cosets contained in them. Let us
briefly recall this here. The point is to describe the cosets contained in $X$. Call
such a coset $b+B$ ($b$ a point, $B$ an abelian subvariety of $A$) {\it maximal}
if it is not included in any other larger torus coset still contained in $X$. Then one
may prove that the set of abelian subvarieties $B$ is finite,  for $b+B$ running
through maximal cosets. A simple proof of this is in [BZ], Lemma 2. (It uses only
rather standard considerations involving degrees.) With this in mind, we proceed
to illustrate the proof of Theorem 1.1* by induction on $\dim A$, the case of
dimension $0$ (or $1$) being indeed trivial. Using exactly the same method of
proof of the weaker version, it suffices, through Theorem 2.1, to deal with the
torsion points which lie in some coset $b+B$ contained in $X$ and having
dimension $>0$. We may then assume that this coset is maximal, and thus that $B$
lies in a certain finite set. Thus for our purposes we may assume that $B$ is
fixed. The quotient $A/B$ has the structure of an abelian variety $A'$ and we
have $\dim A'<\dim A$. Let $\pi:A\to A'$ be the natural projection. The set of
$b\in A$ such that $b+B$ is contained in $X$ is easily seen to be an algebraic
variety, which projects to a variety $X'$ under $\pi$. The coset $b+B$ contains a
torsion point if and only if $\pi(b)$ is torsion on $A'$. Now, the induction
assumption applied to $A',X'$ easily concludes the proof. 
\bigbreak

\noindent {\bf Final remarks. \/} 

1. The Manin-Mumford conjecture for $A$
and $X$ defined over $\CC$ follows from the above by specialization arguments,
as in the original papers of Raynaud. Moreover, it is known that a version that
is uniform (regarding the number of cosets required to contain all the torsion
points) as $X$ varies over a family of subvarieties of fixed dimension and
degree of a given abelian variety follows from the above version, as described
in [BZ] or see Hrushovski [H] or Scanlon [Sc] (so-called ``automatic
uniformity''). 

2. It seems that an argument along the present lines can also be
given for the easier multiplicative version of Manin-Mumford mentioned in \S1.
Here the role of our Theorem 2.1 would be played by the theorem of Ax [Ax]
establishing Schanuel's conjecture for power series. 

3. Quantitative versions of the Manin-Mumford conjecture (and indeed of the Mordell-Lang conjecture) giving upper estimates for the number of torus cosets are given by R\'emond [Re]. Good
dependencies e.g. on the degree of $X$ when $A$ is fixed also follow from Hrushovski's proof [H]. Our method does not for the present yield any new quantitative information. The result of Pila-Wilkie holds in an arbitrary o-minimal structure. In this generality it seems one could not hope for good
dependence e.g. on the degree of $X$ in the Pila-Wilkie result. For the specific
case of sets defined by algebraic relations among theta-functions one might hope
for good bounds, but we have not pursued this. For the multiplicative version,
bounds with good dependencies on the degrees of the varieties were obtained by R\'emond [Re2]; see also Beukers-Smyth (for curves) and Aliev-Smyth (in general) [BS, AS]. 

4. It seems clear that the constant $c_1$ can be taken uniformly as $X$
varies over all subvarieties of $A$ of fixed dimension and degree, by standard
uniformity properties in o-minimal structures.

 5. In a paper in progress by Masser and Zannier, the present method has been applied to prove the following: {\it For $l\neq 0,1$, let $E_l$ be the elliptic curve $y^2=x(x-1)(x-l)$ and let $P_l,Q_l$ be two points on $E_l$ with $x$-coordinate resp. 2,3. Then there are only finitely many complex values of $l$ such that both $P_l,Q_l$ are torsion on $E_l$.}
 This kind of result is related to Silverman's specialization theorem,  a special case of which implies  that the $l\in\CC$ such that $P_l$ or  $Q_l$ is torsion form a set of algebraic numbers of bounded height. The finiteness statement however seems not to follow directly from any known result. (The Manin-Mumford statement, applied to a suitable subvariety of $E\times E$,  would work if $P_l,Q_l$ would be taken as variable  $\ZZ$-independent points on a fixed elliptic curve $E$.)

\bigskip 

\noindent{\bf Acknowledgements.} It is a pleasure to thank David Masser for very helpful discussions and for generously providing us with  the exact  parts of his work most relevant here. We also thank Matt  Baker for his kind  interest and indication of references. The first author is grateful to the Scuola Normale Superiore di Pisa for support and hospitality during the preparation of the present paper, and the second author is similarly grateful to the University of Bristol.
\bigskip

\bigskip

\centerline{\title References}\medskip

\item{[AS]} - I.
Aliev and C. Smyth, Solving equations in roots of unity,
arXiv:0704.1747v3 [math.NT], dated 1 February 2008. 

\item{[Ax]} - J. Ax, On Schanuel's
conjectures. {\it Ann. of Math. (2)\/} {\bf 93\/} (1971), 252--268 

\item{[BR]} - M. Baker and K. Ribet, Galois theory and torsion points on curves. Les XXIImes JournŽes Arithmetiques (Lille, 2001).
{\it J. ThŽor. Nombres Bordeaux\/} {\bf 15\/} (2003),  11--32,

\item{[BS]} -
F. Beukers and C. Smyth, Cyclotomic points on curves, {\it Number theory for the
millenium (Urbana, Illinois, 2000),\/} I, A K Peters, 2002.

 \item {[BM]} - E.
Bierstone and P. Milman, Semianalytic and subanalytic sets, {\it Publ. Math.
IHES\/} {\bf 67\/} (1988), p. 5--42.

 \item{[BP]} - E. Bombieri and J. Pila, The
number of integral points on arcs and ovals, {\it Duke Math. J.}, {\bf 59}
(1989), 337-357.

 \item{[BZ]} - E. Bombieri and U. Zannier, Heights of algebraic
points on subvarieties of abelian varieties, {\it Ann. Scuola Norm. Sup. Pisa
(4)} {\bf 23} (1996), 779--792.

 \item {[G]} - A. Gabrielov, Projections of
semianalytic sets, {\it Funkcional. Anal. i Prilo\v zen.} {\bf 2} (1968),
18--30. English translation: {\it Functional Anal. Appl.} {\bf 2\/} (1968),
282--291. \item{[H]} - E. Hrushovski, The Manin-Mumford conjecture and the model
theory of difference fields, {\it Ann. Pure. Appl. Logic\/} {\bf 112\/} (2001),
43--115.

\item{[L]} - S. Lang, {\it Fundamentals of Diophantine Geometry},
Springer-Verlag, 1983. \item{[M]} - D. Masser, Small values of the quadratic
part of the Neron-Tate height on an abelian variety, {\it Compositio} {\bf 53\/}
(1984), 153--170.

\item{[MZ]} -D. Masser and U. Zannier, Torsion anomalous points and families of elliptic curves, paper in progress.

\item{[O]} - F.
Oort, ``The'' general case of S. Lang's conjecture (after Faltings), {\it
Diophantine approximation on abelian varieties,\/} Edixhoven and Evertse,
editors, LNM 1566, Springer 1993/2003. 

\item{[P]} - J. Pila, Integer points on
the dilation of a subanalytic surface, {\it Quart. J. Math.} {\bf 55} (2004),
207-223. 

\item{[PW]} - J. Pila and A. J. Wilkie, The rational points of a
definable set, {\it Duke Math. J.}, {\bf 33} (2006), 591-616. 

\item{[PR]} - R. Pink and D. Roessler, On Hrushovski's proof of the Manin-Mumford
conjecture, {\it Proceedings of the International Congress of Mathematicians, Vol. I (Beijing, 2002)}, 539--546, Higher Ed. Press, Beijing, 2002.

\item{[Py]} - A.
Pillay, Model Theory and Diophantine Geometry, {\it Bull. A.M.S.}, Vol. 34,
1997, 405-422.

\item{[RU]} - N. Ratazzi and E. Ullmo, Galois + \'equidistribution = Manin-Mumford, web preprint dated
24 September 2007, available at http://www.math.u-psud.fr/$\sim$ullmo/

 \item{[R]} - M. Raynaud, Courbes sur une vari\'et\'e ab\'elienne
et points de torsion, {\it Invent. Math.} { {\bf 71} (1983), no. 1, 207--233.

\item{[R2]} - M. Raynaud, Sous-vari\'et\'es d'une vari\'et\'e ab\'elienne et
points de torsion, in {\it Arithmetic and Geometry}, Vol. I, Birkhauser, 1983.

\item{[Re]} - G. R\'emond, Decompte dans un conjecture de Lang, {\it Invent.
Math.} {\bf 142} (2000), 513--545. 

\item{[Re2]} - G. R\'emond, Sur les sous-vari\'et\'es des tores, {\it Compsotio\/} {\bf 134\/} (2002), 337--366.

\item{[Ro]} - M. Rosen, Abelian varieties
over $\CC$, {\it Arithmetic geometry,} Cornell and Silverman, editors, Springer,
1986.

\item{[Sc]} - T. Scanlon, Automatic uniformity, {\it International Mathematics
Research Notices\/} 2004, no. 62, 3317 - 3326. 

\item{[Sh]} - M. Shiota, {\it
Geometry of Subanalytic and Semialgebraic Sets}, Birkh\"auser, 1997. 

\item{[S]} - C.-L. Siegel, {\it
Lectures on the geometry of numbers,} Springer, Berlin, 1989.

 \item{[T]} - P.
Tzermias, The Manin-Mumford conjecture: a brief survey, {\it Bull. London Math.
Soc.\/} {\bf 32\/} (2000), 641--652.

 \bigskip \line{School of Mathematics \hfill
Scuola Normale Superiore}
\line{University of Bristol \hfil Piazza dei Cavalieri 7}
 \line{Bristol BS8 1TW \hfill 56126 Pisa} \line{UK \hfill Italy}
\bye